\documentclass[12pt]{amsart}
\usepackage{amssymb,amsmath} 
\newif\ifmypdf
\mypdffalse
\ifmypdf\usepackage[pdftex]{graphicx}\else \usepackage[dvips]{graphicx}\fi

\newcommand{\Graphics}[2]{\ifmypdf \includegraphics[#1]{#2.pdf} \else \includegraphics[#1]{#2.eps}\fi}

      \newenvironment{changemargin}[2]{\begin{list}{}{
         \setlength{\topsep}{0pt}\setlength{\leftmargin}{0pt}
         \setlength{\rightmargin}{0pt}
         \setlength{\listparindent}{\parindent}
         \setlength{\itemindent}{\parindent}
         \setlength{\parsep}{0pt plus 1pt}
         \addtolength{\leftmargin}{#1}\addtolength{\rightmargin}{#2}
         }\item }{\end{list}}

\newcommand{\LE}{\preccurlyeq}
\newcommand{\MajLE}{\LE_{_\mathrm{Maj}}}
\newcommand{\AvLE}{\LE_{_\mathrm{Av}}}

\newcommand{\myfigure}[1]{\goodbreak\begin{figure}[!htp]#1\end{figure}}
\newcommand{\mytable}[1]{\goodbreak\begin{table}[!htp]\begin{changemargin}{-2cm}{-2cm}\begin{center}
#1\end{center}\end{changemargin}\end{table}}
\newtheorem{thm}{Theorem}[section]
\newtheorem{prop}[thm]{Proposition}
\newtheorem{cor}[thm]{Corollary}
\newtheorem{lem}[thm]{Lemma}
\newtheorem{conj}[thm]{Conjecture}
\newtheorem{exa}[thm]{Example}

\newtheorem{notat}[thm]{Notation}
\newtheorem{alg}[thm]{Algorithm}
\newtheorem{prob}[thm]{Problem}
\theoremstyle{definition}
\newtheorem{defn}[thm]{Definition}
\theoremstyle{remark}
\newtheorem{rem}[thm]{Remark}

\newcommand{\inv}{^{-1}}

\newcommand{\be}{\begin{enumerate}}
\newcommand{\ee}{\end{enumerate}}
\renewcommand{\i}{\item}
\newcommand{\ble}{\begin{lem}}
\newcommand{\ele}{\end{lem}}
\newcommand{\bth}{\begin{thm}}
\newcommand{\bpr}{\begin{prop}}
\newcommand{\epr}{\end{prop}}
\newcommand{\bco}{\begin{cor}}
\newcommand{\eco}{\end{cor}}
\newcommand{\bcon}{\begin{conj}}
\newcommand{\econ}{\end{conj}}
\newcommand{\bde}{\begin{defn}}
\newcommand{\ede}{\end{defn}}
\newcommand{\bex}{\begin{exa}}
\newcommand{\eex}{\end{exa}}
\newcommand{\brem}{\begin{rem}}
\newcommand{\erem}{\end{rem}}
\newcommand{\bnot}{\begin{notat}}
\newcommand{\enot}{\end{notat}}
\newcommand{\balg}{\begin{alg}}
\newcommand{\ealg}{\end{alg}}

\newcommand{\<}{\left <}
\renewcommand{\>}{\right > }

\newcommand{\si}{\sigma}

\newcommand{\cB}{\mathcal{B}}

\long\def\forget#1\forgotten{} %
\newcommand{\Garlen}{{\ell_\mathrm{G}}}
\newcommand{\RedGarlen}{{\ell_\mathrm{RG}}}

\begin{document}

\title[Length-based conjugacy search in the braid group]{Length-based conjugacy search in the braid group}

\author[Garber, Kaplan, Teicher, Tsaban, Vishne]{David Garber,
  Shmuel Kaplan, Mina Teicher, Boaz Tsaban, and Uzi Vishne}

\address{David Garber,
Einstein institute of Mathematics, The Hebrew University,
Givat-Ram 91904, Jerusalem, Israel; and Department of Sciences,
Holon Academic Institute of Technology, 52 Golomb Street, Holon
58102, Israel} \email{garber@math.huji.ac.il, garber@hait.ac.il}

\address{Shmuel Kaplan, Mina Teicher, and Uzi Vishne,
Department of Mathematics and Statistics, Bar-Ilan University,
Ramat-Gan 52900, Israel}
\email{[kaplansh, teicher, vishne]@math.biu.ac.il}

\address{Boaz Tsaban, Department of Mathematics,
Weizmann Institute of Science, Rehovot 76100, Israel}
\email{boaz.tsaban@weizmann.ac.il}
\urladdr{http://www.cs.biu.ac.il/\~{}tsaban}

\thanks{\tiny This paper is a part of the Ph.D.\ thesis of the second named author
 at Bar-Ilan University.}
\thanks{\tiny
This research was partially supported by the Israel Science
Foundation through an equipment grant to the school of Computer
Science in Tel-Aviv University. The authors were partially
supported by: Golda Meir Fellowship (first named author),
EU-network HPRN-CT-2009-00099(EAGER), ENRI, the Minerva
Foundation, and ISF grant \#8008/02-3 (second and third named
authors), the Koshland Center for Basic Research (fourth named
author).}

\begin{abstract}
Several key agreement protocols are based on the following
\emph{Generalized Conjugacy Search Problem}: Find,
given elements
$$b_1,\dots,b_n\mbox{ and }xb_1x\inv,\dots,xb_nx\inv$$
in a nonabelian group $G$, the conjugator $x$. In the case of
subgroups of the braid group $B_N$, Hughes and Tannenbaum
suggested a length-based approach to finding $x$. Since the
introduction of this approach, its effectiveness and
successfulness were debated.

We introduce several effective realizations of
this approach.
In particular, a length function is defined on
$B_N$ which possesses significantly better properties
than the natural length associated to the Garside
normal form. We give experimental results concerning
the success probability of this approach, which suggest
that an unfeasible computational power is required
for this method to successfully solve the Generalized Conjugacy Search Problem
when its parameters are as in existing protocols.
\end{abstract}

\maketitle

\section{Introduction}
Assume that $G$ is a nonabelian group.
The following problem has a long history and many applications
(see \cite{BGC}).
\begin{prob}[Generalized Conjugacy Search Problem]\label{GCSP}
Given elements $b_1,\dots,b_n\in G$
and their conjugations by an unknown element $x\in G$,
$$xb_1x\inv,xb_2x\inv,\dots,xb_nx\inv,$$
find $x$ (or any element $\tilde x\in G$ such that
$\tilde xb_i\tilde x\inv = xb_ix\inv$ for $i=1,\dots,n$.)
\end{prob}
In the sequel, we will not make any distinction between the actual conjugator
$x$ and any other conjugator $\tilde x$ yielding the same results.

The \emph{braid group} $B_N$ is the group generated by the
$N-1$ \emph{Artin generators} $\si_1,\dots,\si_{N-1}$,
with the relations
\begin{eqnarray*}
\si_i \si_{i+1}\si_i & = & \si_{i+1} \si_i \si_{i+1},\\
\si_i \si_j          & = & \si_j \si_i \mbox{ when }|i-j|>1
\end{eqnarray*}
Information on the basic algorithms in the braid group is
available in \cite{Cha} and the references therein. We will focus
on the case where $G$ is the subgroup of $B_N$ generated by given
elements $a_1,\dots,a_m$. A solution of the generalized conjugacy
problem in this case immediately implies the vulnerability of
several cryptosystems introduced in \cite{Anshel, Ko},
and the methods of solution may be applicable
to several other cryptosystems from \cite{Anshel, Paeng}.

\subsection*{History, motivation, and related work}
The length-based approach to the Conjugacy Problem was suggested by Hughes and Tannenbaum in \cite{HT},
as a potential attack on the cryptosystems introduced in \cite{Anshel, Ko}.
Based on \cite{HT}, Garrett \cite{Garrett} has doubted
the security of these cryptosystems.
But soon afterwards he published an errata withdrawing
these doubts (see \cite{BGC}). The reason was that no known realization of
Hughes and Tannenbaum's scheme
(i.e., definition of actual, effective length functions) was given before,
and in particular, the success probability of this approach could not be
estimated.
The purpose of the current paper is to introduce and compare several
such realizations, and provide actual success probabilities for specific parameters.

We stress that we are not interested here in the best possible solution
of the generalized conjugacy problem, but rather in settling the
debate concerning the applicability of the Hughes-Tannenbaum length-based
approach to the problem.

Other approaches appear in \cite{HofSte, LeePark} and turn out more successful.
However, the length-based approach has several advantages:
First, one does not need to know the conjugated element in order
to find the conjugator using this approach, and second,
it essentially deals with \emph{arbitrary equations}.
The current paper gives the foundations of this approach,
on which we build in \cite{BraidEqns}, where an extension of this
approach is suggested and good success rates are achieved for arbitrary
equations.

Some of the citations of the present paper (see \cite{Bangert, DeBBC, Gon, Lee, Shp, UshShp, ShpZap, DehShift})
refer to its preliminary draft \cite{draft}, which contains much more details
and examples. We have tried to make the present version concise.

\subsection*{Length-based attacks}\label{LengthBasedAttacks}
Throughout this paper we make the following assumptions:
\be
\i The conjugator $x$ belongs to
a given finitely generated subgroup of $B_N$, whose generators
$$\{a_1,\dots,a_m,a_1\inv,\dots,a_m\inv\}$$
are given,
\i $x$ was generated as a product of a fixed, known number of
generators $a_i^{\pm 1}$ chosen at random from the set of
generators;
\i We are given elements
$$xb_1x\inv,\dots,xb_nx\inv$$
where each $b_i\in B_N$ is generated by some (nontrivial) random process,
and we wish to find $x$.
\ee
We try to find the conjugator $x$
by using the property that for an appropriate, efficiently computable
\emph{length function} $\ell$ defined on $B_N$,
$\ell(a\inv ba)$ is usually greater than $\ell(b)$ for elements $a,b\in B_N$.
Therefore, we try to reveal $x$ by
peeling off generator after generator from the given braid elements
$xb_1x\inv,\dots,xb_nx\inv$:
Assume that
$$x = g_1\cdot g_2\cdots g_k,$$
where each $g_i$ is a \emph{generator}. We fix some
\emph{linear order} $\LE$ on the set of all possible $n$-tuples of
lengths, and choose a generator $g$ for which the lengths vector
$$\<\ell(g\inv xb_1x\inv g), \dots, \ell(g\inv xb_nx\inv g)\>$$
is minimal with respect to $\LE$.
With some nontrivial probability,
$g$ is equal to $g_1$ (or at least, $x$ can be rewritten as a product of $k$ or fewer
generators such that $g$ is the first generator in this product),
so that $g\inv x = g_2\cdots g_k$ is a product of fewer generators and we may continue this way,
until we get all $g_i$'s forming $x$.

If one is capable of doing
$O((2m)^t)$ computations, it is better to check all
possibilities of $g_1\cdots g_t$ by peeling off $g_1\cdots g_t$
from $x$ and choosing the $t$-tuple which yielded the minimal
lengths vector. We will call this approach \emph{look ahead of depth $t$}.

In order for any of the above to be meaningful, we must define the
length function $\ell$ and the linear ordering $\LE$.
We will consider several candidates for these.

\section{Realizations of the length function}

We assume that each generator
$a_i$ is obtained by taking a product of some fixed
number of (randomly chosen) Artin generators, to whom we refer as the
``length'' of the generators.
Unless otherwise stated, in all of our experiments
the length of each element $a_i$ is $10$.
By a \emph{generator} we mean either an element $a_i$, $i=1,\dots,m$,
or its inverse.
We will (informally) write $|x|=n$ when we mean that $x$ was generated by a product
of $n$ generators chosen at random (with uniform distribution)
from the list of $2m$ generators $a_1,\dots,a_m,a_1\inv,\dots,a_m\inv$.

\subsection{The length function $\ell$} \label{LEN_FUNC}
The \emph{Garside normal form} of an element $w\in B_N$
is the unique presentation of $w$ in the form
$\Delta_N^{-r} \cdot p_1\cdots p_k$,
where $r \geq 0$ is minimal and $p_1,\dots,p_k$ are
permutation braids in left canonical form \cite{Cha}.
Using the Garside normal form, one can assign a ``length'' to each
$w\in B_N$ efficiently \cite{Cha}.

\bde
The \emph{Garside length} of an element $w \in B_N$, $\Garlen(w)$,
is the number of Artin generators needed to write $w$ in its
Garside normal form.
If the Garside normal form of $w$ is $\Delta_N^{-r} \cdot p_1\cdots p_k$, then
$$\Garlen(w)=r\cdot \binom{N}{2}+\sum_{i=1}^{k}{|p_i|},$$
where $|p|$ denotes the length of the permutation $p$.\footnote{The
length of a permutation $p$ is the number
of order distortions in $p$, that is, pairs $(i,j)$ such that $i<j$ and $p(i)>p(j)$.
}
\ede

The problem with this function is that it is not close enough to being monotone with $|x|$:
One has to multiply many generators before an increase
in the length function is observed.
The left part of Figure \ref{GLgraph} shows, for a fixed word $b$,
$\Garlen(xbx\inv)$ as a function of $|x|$.
Its right part shows the \emph{average} of
$\Garlen(xbx\inv)$ computed over $1200$ random words.

\myfigure{
\begin{changemargin}{-2cm}{-2cm}
\begin{center}
\Graphics{scale=0.58}{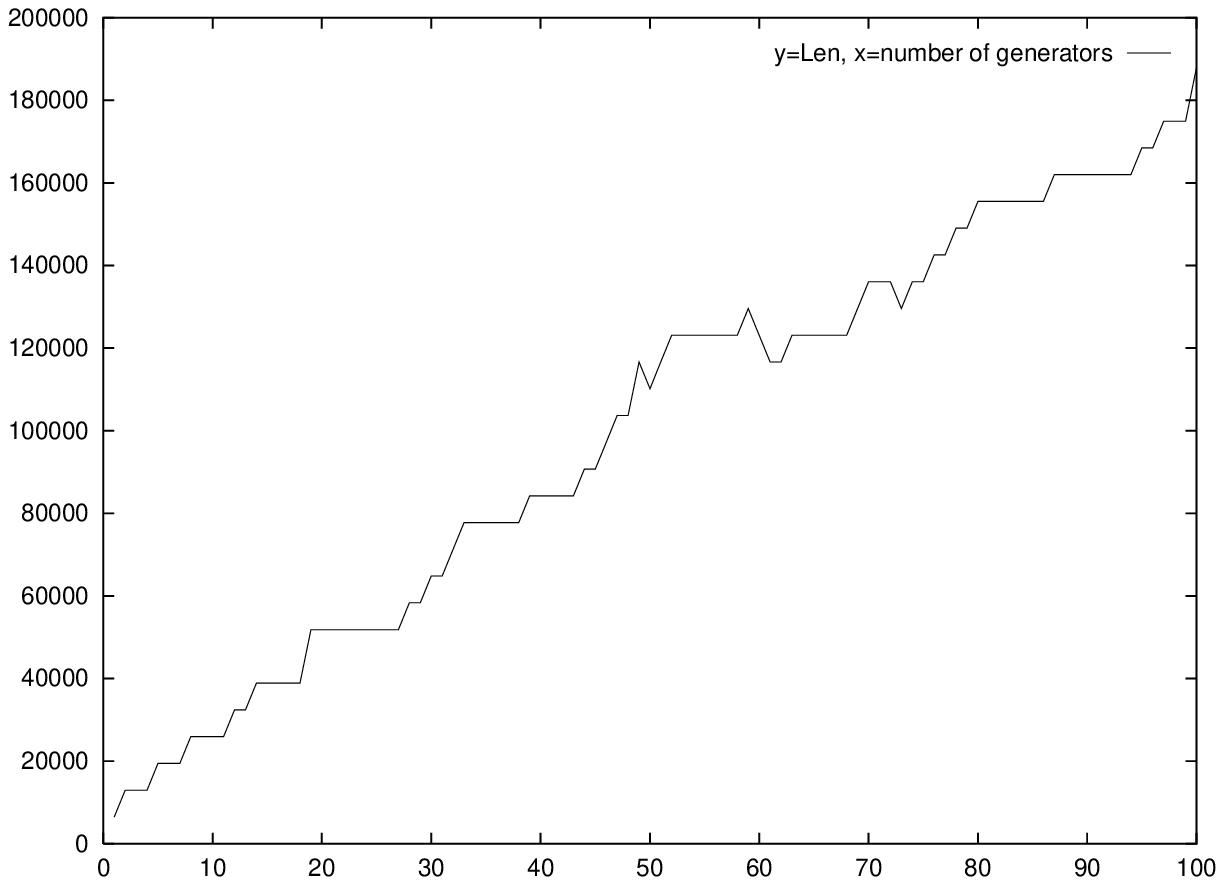}
\Graphics{scale=0.58}{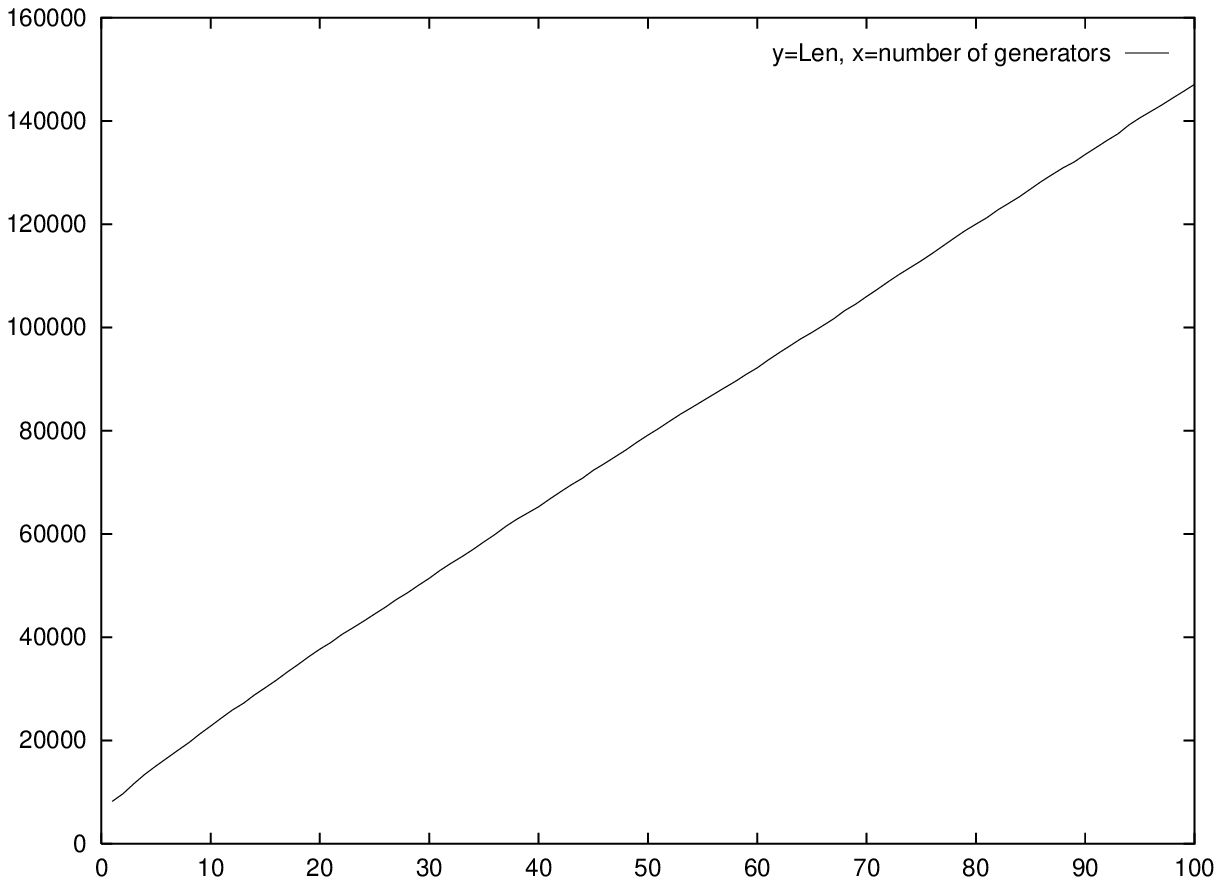}\\
\caption{The growth of $\Garlen(w)$: Specific case (left) and average growth (right)}\label{GLgraph}
\end{center}
\end{changemargin}
}

We wish to have a length function that is closer to being monotone.
For each permutation braid $p$,
$\tilde p:=p\inv\Delta_N$ is a permutation braid.
Thus, if $w=\Delta_N^{-r} \cdot p_1 \cdots p_k$
and $r>0$, we can replace $\Delta_N\inv p_1$ with
$\tilde p_1\inv$ to get
$w=\Delta_N^{-(r-1)} \cdot \tilde p_1\inv p_2 \cdots p_k$.
Now, $\Delta_N$ almost commutes with any permutation
braid: For each permutation braid $q$ there exists a permutation braid $q'$
such that $|q'|=|q|$ and $q\Delta_N = \Delta_N q'$, that is,
$\Delta_N\inv q\inv = (q')\inv\Delta_N\inv$.
Consequently,
$w=\Delta_N^{-(r-2)} \cdot (\tilde p_1')\inv \Delta_N\inv p_2 \cdots p_k$,
and we can replace $\Delta_N\inv p_2$ with $\tilde p_2\inv$ as before.
We iterate this process as much as possible, to get a presentation
$$w = \begin{cases}
\Delta_N^{-(r-k)} (\tilde p_1')\inv\cdots (\tilde p_k')\inv & k<r\\
(\tilde p_1')\inv\cdots (\tilde p_r')\inv\cdot p_{r+1}\cdots p_k  & r\le k
\end{cases}$$
In each case, $w$ has the form $a\inv b$ where $a,b$ are positive braid
words or the identity element, and we define the \emph{reduced Garside length}
to be the sum of the length of $a$ and the length of $b$.\footnote{The length of a positive braid
word is well defined to be the number of generators in its presentation.}
This is equivalent to the following.
\bde
Let $w=\Delta_N^{-r} \cdot p_1 \cdots p_k$ be the Garside normal
form of $w$. The \emph{Reduced Garside length} of $w$  is defined by
$$\RedGarlen(w) = \Garlen(w)-2 \sum_{i=1}^{\min(r,k)}{|p_i|}$$
\ede
This function turns out much closer to monotone than $\Garlen$ --
see Figure \ref{RGLgraph}.

\myfigure{
\begin{changemargin}{-2cm}{-2cm}
\begin{center}
\Graphics{scale=0.58}{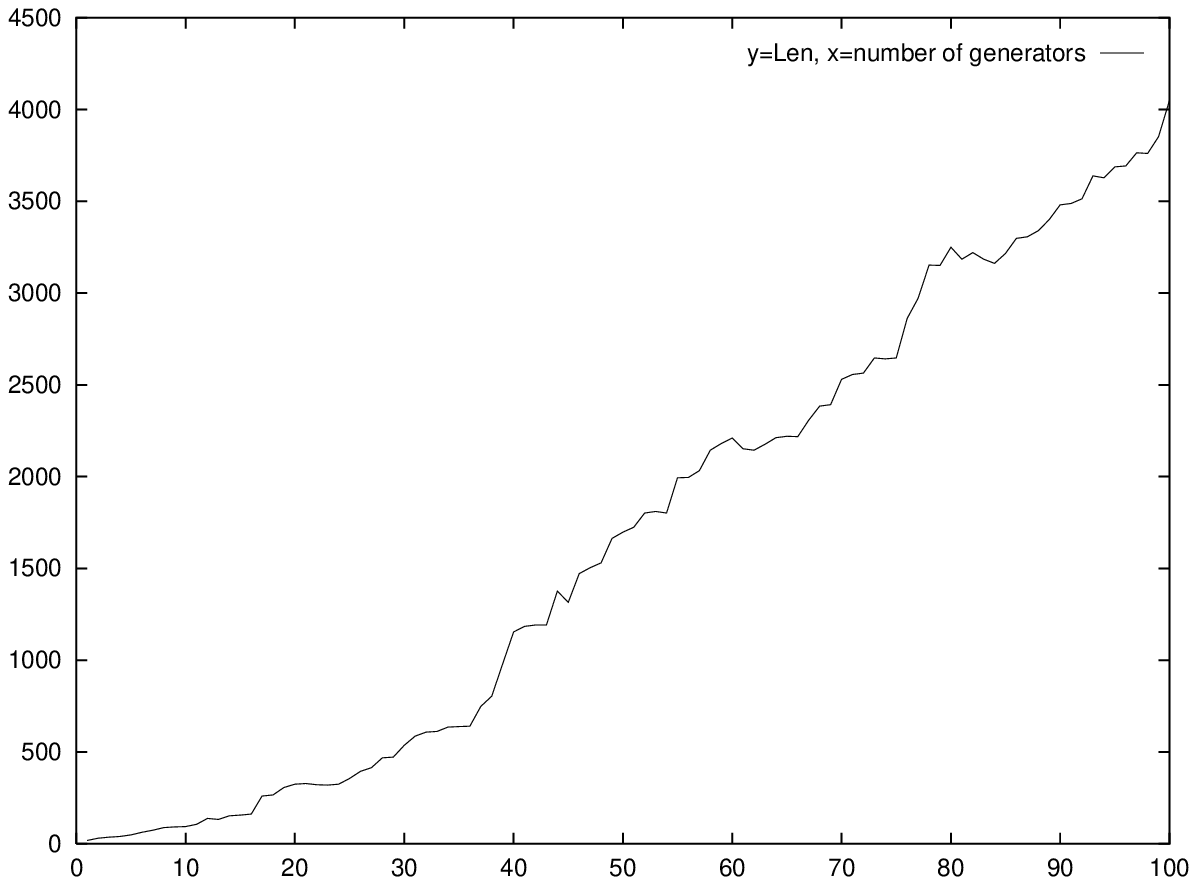}
\Graphics{scale=0.58}{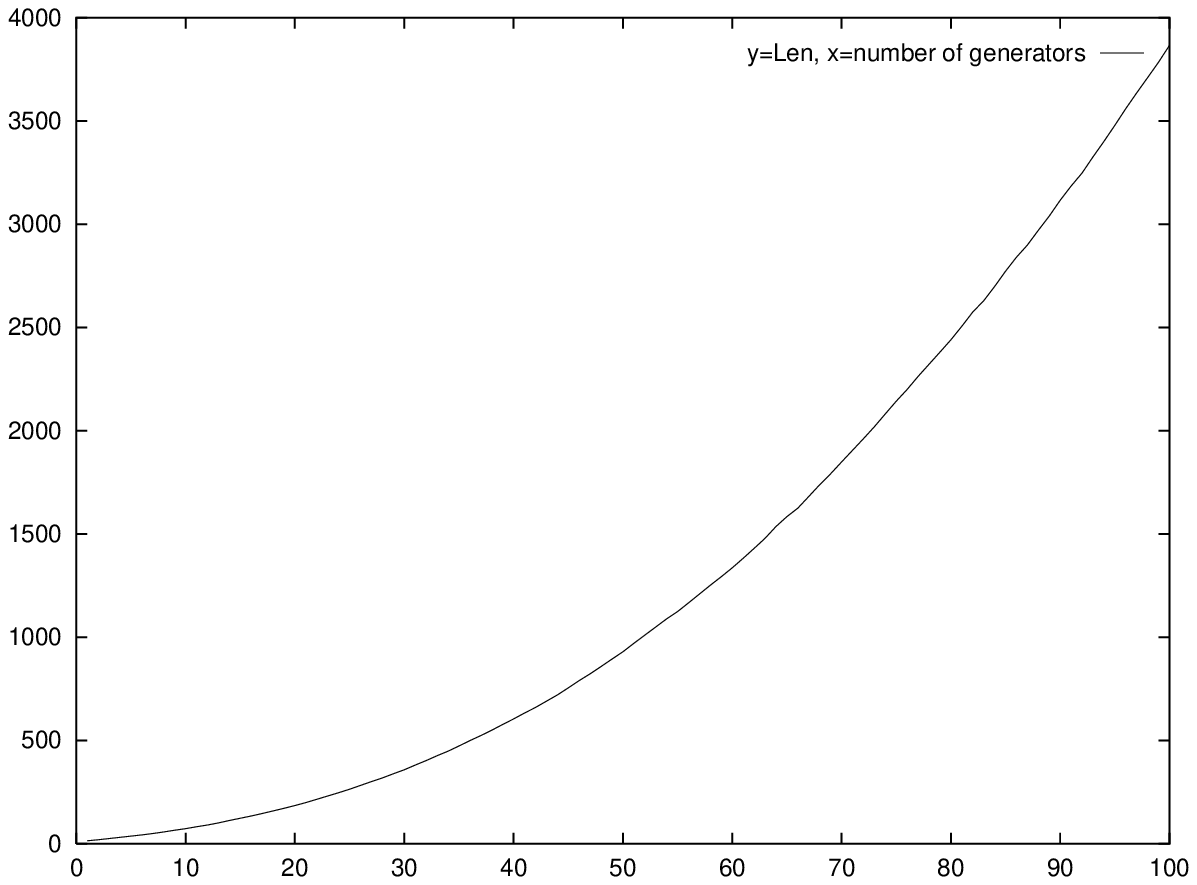}\\
\caption{The growth of $\RedGarlen(w)$: Specific case (left) and average growth (right)}\label{RGLgraph}
\end{center}
\end{changemargin}
}

\subsection{Statistical comparison of the length functions}
The purpose of the length function $\ell$ is to distinguish
between the case $|X|=k-1$ (after peeling off a correct generator)
and $|X|=k+1$ (after trying to peel off a wrong generator).
Thus, a natural measure for the effectiveness of the length function is
the distance in standard deviations between $\ell(X')$ and $\ell(X)$
when $|X'|=|X|+2$.

We fixed a random set of $20$ generators in $B_{81}$, and computed
(an approximation of)
$E(\ell(X')-\ell(X))/\sqrt{V(\ell(X')-\ell(X))}$ as a function of
$|X|$ for $|X|=1,\dots,100$.
(Roughly speaking, when $n$ independent samples are added,
the effectiveness of the comparison is $\sqrt{n}$ times this number.)
We did that for both $\Garlen$ and
$\RedGarlen$. The results appear in Figure
\ref{E/sqrt(V) sp sbgp}, and show that the score for $\RedGarlen$ is
significantly higher. This phenomenon is typical -- we have checked several random
subgroups of the braid group and all of them exhibited the same behavior.
\myfigure{
\begin{changemargin}{-2cm}{-2cm}
\begin{center}
\Graphics{scale=0.38}{graph7}
\caption{Distance between right and wrong in standard deviations.}\label{E/sqrt(V) sp sbgp}
\end{center}
\end{changemargin}
}

More evidence for the superiority of $\RedGarlen$ over $\Garlen$ will be given in the
following sections.

\section{Realizations of the linear ordering $\LE$}\label{linearorder}
Recall that after peeling off a candidate for a generator and evaluating
the resulting lengths, we need to compare the vectors of lengths
according to some linear ordering $\LE$,
and choose a generator which gave a minimal vector with respect
to $\LE$.
We tested two natural linear orderings.

The most natural approach is to take the average of the lengths in the vector.
This is equivalent to the following.
\begin{defn}[Average based linear ordering]
$$\<\alpha_1,\dots,\alpha_n\>\AvLE\<\beta_1,\dots,\beta_n\>\quad\mbox{if}\qquad
\sum_{i=1}^n\alpha_i\le \sum_{i=1}^n\beta_i.$$
\end{defn}

With this at hand, we have performed the following experiment.
We fixed a subgroup of $B_{81}$ generated by $m=20$ generators.
Then we chose at random $200$ elements of the form $xw_j$ which share
the same leading prefix $x$, and for each generator $a_i^{\pm 1}$ we computed
$\ell(a_i^{\pm 1}xw_j)$ for each $j$ (and $\ell=\Garlen$ or $\RedGarlen$).
For each of these two length functions $\ell$,
we have sorted the resulting length vectors according to $\AvLE$
and checked the position of the ``correct'' generator, i.e., the
generator which appeared leftmost in our computation of the word.\footnote{In principle
there could be more than one ``correct'' generator, but when the generators
are long enough this is unlikely to happen often.}
We repeated the computations for $138$ distinct $X$'s,
and for $|X|=40$ and $|X|=100$.
For an ideal length function (and an ideal linear ordering $\LE$),
the correct generator
would always be ranked first, and the results in Figure \ref{nicegraph}
show that $\RedGarlen$ is closer to this ideal than $\Garlen$:
In the graphs, we show the \emph{distribution} (lower part of the graph) and the
\emph{accumulated distribution} (upper part of the graph) of the position of the
correct generator, for each of the length functions.

\myfigure{
\begin{changemargin}{-2cm}{-2cm}
\begin{center}
\Graphics{scale=0.37}{Len40all}
\Graphics{scale=0.37}{Len100all}\\
\caption{Position of correct generator $\Garlen$ and
$\RedGarlen$}\label{nicegraph}
\end{center}
\end{changemargin}
}

However, it turns out that even for the better length function $\RedGarlen$,
the task of identifying the correct generator is not trivial.
To demonstrate this, we selected at random one of the cases of $x$ from the
previous experiment, and computed over the given $200$
samples the distribution of $\RedGarlen$ for
each generator. Figure \ref{realdistrib} shows the distribution
for the correct generator (in boldface) and of arbitrarily chosen
$7$ out of the remaining $40$ generators (for an aesthetic reason we did not
plot all $40$).

\myfigure{
\begin{changemargin}{-2cm}{-2cm}
\begin{center}
\Graphics{scale=0.35}{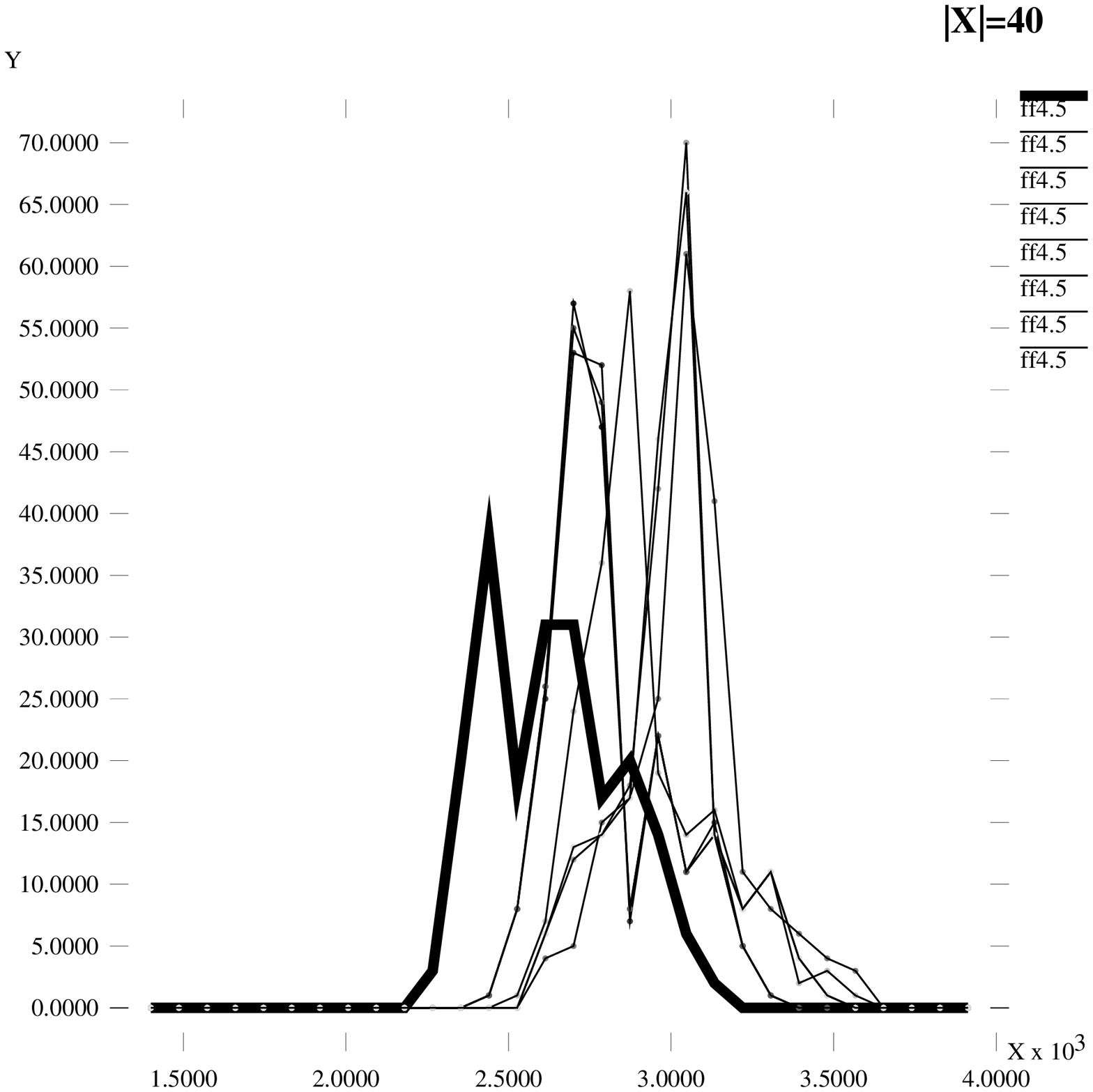}
\Graphics{scale=0.35}{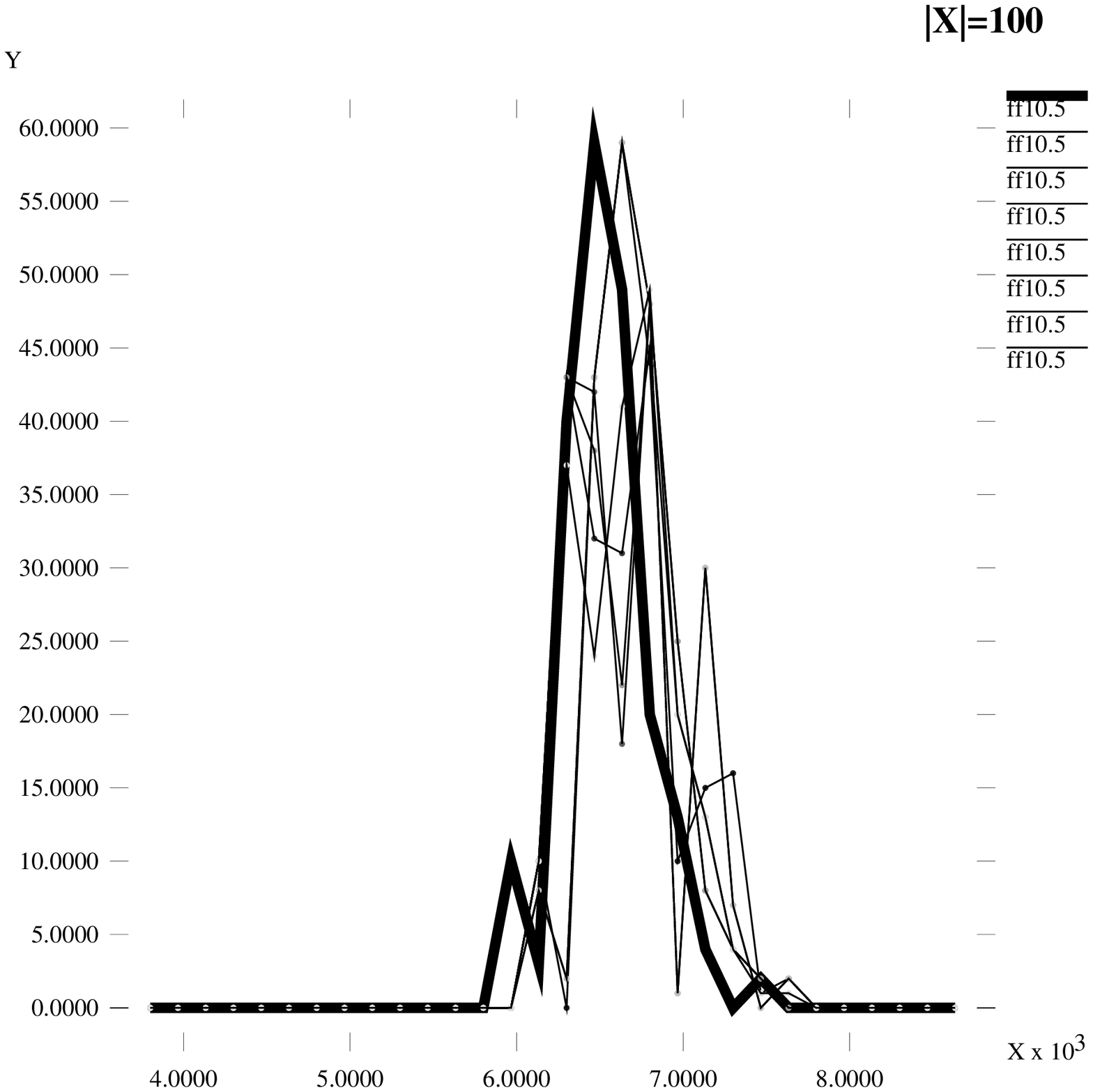}\\
\caption{Actual distribution.}\label{realdistrib}
\end{center}
\end{changemargin}
}

While the correct distribution tends more to the left (i.e., to smaller values),
there is a large overlap with the rest of the distributions.
We must emphasize that while Figure \ref{realdistrib} demonstrates the
typical case, there exist cases where the distribution of the correct
generator is not the leftmost. In these cases the current method is doomed to
fail, no matter how many conjugations we are given for the same conjugator.

Finally, for the sake of comparison, we define one more natural linear ordering
of the space of length vectors.
We expect the correct generator to yield the shortest length more
often than the other generators. This motivates the following definition.
\begin{defn}[Majority based linear ordering]
Consider the set of all obtained length vectors.
For each $i=1,\dots,n$, consider the $i$th coordinate of each vector and
let $\mu_i$ denote the minimum of all these $i$th coordinate values.
Then
$$\<\alpha_1,\dots,\alpha_n\>\MajLE\<\beta_1,\dots,\beta_n\>\quad\mbox{if}\qquad
|\{i : \alpha_i = \mu_i\}| \ge |\{i: \beta_i=\mu_i\}|.$$
\end{defn}
In the following section we compare the success probabilities of the length-based
approach using the two length functions and two linear orderings defined in this
section.

\section{Experimental results for the conjugacy problem}

\subsection{The probability of obtaining the correct generator}\label{find_correct_prob}
In this experiment we determine the probability that the correct generator is indeed
the minimal with respect to the length function $\ell$ and linear ordering $\LE$ used.
The choice of parameters in the experiments throughout the paper are usually motivated by the choices
given in \cite{Anshel}, which are believed there to make the generalized conjugacy
problem difficult.

We made $200$ experiments using the following parameters:
$N=81$, $n$ and $m$ (the number of $a_i$'s and $b_i$'s, respectively) are both $20$,
the elements $a_i$ and $b_i$ are products of $10$ random Artin generators,
and $x$ is the product of $5$, $10$, $20$, $40$, $60$, or  $100$
random generators $a_i^{\pm 1}$, respectively.
We tested look ahead depth $t=1,2$.
In each cell of Table \ref{table_B81}, below the
probability that the correct generator is first, we wrote the
probability of its being second.

\mytable{
\begin{tabular}{|l||l|l|l|l|l|l|}
\hline                      & 5      &   10   &   20   &   40   &   60    &   100   \\
\hline
\hline
       $\Garlen$, $\AvLE$, $t=1$           & 0.56   & 0.478  & 0.322  & 0.267  & 0.233   & 0.156   \\
             & 0.16   & 0.188  & 0.1    & 0.167  & 0.089    & 0.1     \\
\hline
       $\Garlen$, $\MajLE$, $t=1$  & 0.43   & 0.344  & 0.222  & 0.244  & 0.178    & 0.156   \\
             & 0.14   & 0.178  & 0.144  & 0.122  & 0.1     & 0.044   \\
\hline
\hline
       $\RedGarlen$, $\AvLE$, $t=1$      & 0.74   & 0.589  & 0.567  & 0.456  & 0.311   & 0.233   \\
        & 0.13   & 0.233  & 0.189  & 0.122  & 0.167   & 0.167   \\
\hline
       $\RedGarlen$, $\MajLE$, $t=1$      & 0.71   & 0.578  & 0.578  & 0.433  & 0.289   & 0.211   \\
         & 0.15   & 0.267  & 0.133  & 0.089  & 0.167   & 0.167   \\
\hline
\hline
       $\Garlen$, $\AvLE$, $t=2$   & 0.433  & 0.287  & 0.111  & 0.1    & 0.114  & 0.099   \\
          & 0.156  & 0.08   & 0.087  & 0.038  & 0.055  & 0.035   \\
\hline
       $\Garlen$, $\MajLE$, $t=2$          & 0.25   & 0.147  & 0.103  & 0.058  & 0.086  & 0.03    \\
              & 0.033  & 0.036  & 0.024  & 0.008  & 0.023  & 0.03    \\
\hline
\hline
       $\RedGarlen$, $\AvLE$, $t=2$      & 0.578  & 0.526  & 0.333  & 0.242  & 0.2    & 0.168   \\
        & 0.183  & 0.127  & 0.135  & 0.138  & 0.105  & 0.05    \\
\hline
       $\RedGarlen$, $\MajLE$, $t=2$      & 0.511  & 0.482  & 0.31   & 0.242  & 0.186  & 0.149   \\
        & 0.139  & 0.139  & 0.127  & 0.104  & 0.091  & 0.089   \\
\hline
\end{tabular}
\caption{The probability that the correct generator is first or second}\label{table_B81}
}

Table \ref{table_B81} shows
that the Reduced Garside length function $\RedGarlen$ is significantly better than the standard
Garside length function $\Garlen$.
Also, observe that using look ahead depth $2$ is preferable to using look ahead depth $1$
twice (to see this, square the probabilities for $t=1$).
Another natural approach to using look ahead $t>1$ is to consider only the first generator
(of the word with the least score) as correct, and ignore the rest of the generators.
This means that in the algorithm for finding $x$, we peel off only one generator at a time
despite the fact that we used look ahead $t>1$. This gives better success rates than
just taking $t=1$, and our experiments indicate that this approach may be slightly better than
that of taking the whole look ahead word, but we did not extensively check this conjecture
since the differences were not significant.
Some other variants of the usage of look ahead are mentioned in \cite{draft}.

\subsection{Nonsymmetric parameters}\label{LONGSHORT}
This experiment checks the effect on the probability of
success when the lengths of the generators $a_i$ and elements
$b_i$ (in terms of Artin generators)
are not equal.

We tested the probability of success for $N=81$, $n=m=20$, look ahead depth $t=2$,
and $|x|=30$.

\begin{table}[!ht]
\centering
\begin{tabular}{|c||c||c|c|c|c|c||c|c|c|c|c|}
\hline
           &         & \multicolumn{5}{c||}{$\RedGarlen$} & \multicolumn{5}{c|}{$\Garlen$} \\
\cline{3-12}
           & length of $b_i$ & \multicolumn{5}{c||}{$a_i$ of length:}     &\multicolumn{5}{c|}{$a_i$ of length:}     \\
\cline{3-12}
           &         &  5  &  10 &  15 &  20 &  25      &  5  &  10 &  15 &  20 &  25     \\
\hline
\hline
$\AvLE$    &    5    &  44 &  82 & 124 & 134 &  156     &  32 &  51 &  81 & 102 & 115     \\
\cline{2-12}
           &   10    &  59 &  97 & 113 & 141 &  150     &  56 &  69 &  79 &  91 &  96     \\
\cline{2-12}
           &   15    &  56 &  91 & 123 & 136 &  141     &  31 &  53 &  75 &  93 &  105    \\
\cline{2-12}
           &   20    &  49 &  77 & 115 & 132 &  149     &  31 &  49 &  77 &  86 &  107    \\
\cline{2-12}
           &   25    &  56 &  84 & 102 & 127 &  141     &  42 &  59 &  60 &  91 &  100    \\
\hline
\hline
$\MajLE$   &    5    &  39 &  70 & 121 & 134 &  160     &  28 &  41 &  66 &  84 &  87     \\
\cline{2-12}
           &   10    &  57 &  97 & 114 & 140 &  156     &  49 &  49 &  58 &  83 &  82     \\
\cline{2-12}
           &   15    &  50 &  85 & 118 & 136 &  144     &  19 &  45 &  59 &  73 &  89     \\
\cline{2-12}
           &   20    &  48 &  80 & 116 & 133 &  149     &  39 &  41 &  52 &  73 &  86     \\
\cline{2-12}
           &   25    &  60 &  89 & 101 & 141 &  152     &  39 &  50 &  56 &  72 &  78     \\
\hline
\end{tabular}
\caption{Number of success out of $200$ tries for different lengths}\label{longshort}
\end{table}

As expected, Table \ref{longshort} shows that if the length of the elements
$a_i$ increases then so does the probability to find a correct generator
(this is like making look ahead deeper without exponentially increasing the number of
candidates for the prefix of $x$). On the other hand, the effect of the length of
the elements $b_i$ is not significant.

\subsection{Increasing the number of given conjugates}\label{many}
Several experiments showed that increasing the number $n$ of given
elements $xb_ix\inv$ from few (about $10$) to many (about $3000$)
did not significantly increase the probability that the correct
generator appears first.

In an instance of the problem the length function $\ell$ and
the (unknown) element $x$ are fixed, and this defines
for each generator $g$ the distribution $F_g$ of
$\ell(g\inv xb x\inv g)$ over random words $b$ of a fixed given length (in
terms of Artin generators).
For each $g$, we have a
sample of the distribution $F_g$ for each given equation.
In most cases, the expectancy of $F_g$ where $g$ is the first letter
in $x$ is smaller than the other expectancies (see Section
\ref{linearorder}), and then enough samples will allow us to
identify $g$. However in some cases the minimal expectancy is
obtained for another generator. In these cases adding more samples
cannot help, and so the probability to find the correct generator
does not tend to $1$ when we increase the number of samples.

Another important observation is that few samples (about
$15$) are needed in order to get very close to the expectancy
of the distributions $F_g$.
In light of the preceding paragraph,
the outcome of the algorithm can be decided
after a relatively small number of samples (i.e., given conjugates)
are collected.
In particular, the success probability
does not significantly improve when $n$ is large.

\subsection{Finding $x$} \label{FIND_x}
The simplest way to try and obtain all generators of $x$ and therefore $x$ would be to
use any of the above algorithms iteratively, at each step peeling off the first
generator.
In the following experiment, the probability to find \emph{all of} $x$ this way was tested.
Here too, the lengths of the $a_i$'s and $b_i$'s were $10$ Artin generators.
We made $500$ experiments, using a weaker variant of $\RedGarlen$ as the length function,
and with no look ahead ($t=1$).
We repeated this for $B_4,\cdots,B_{20}$ and $x$ of lengths $2$ to $18$ generators $a_i^{\pm 1}$.
The result is the number of successes out of $500$ tries.

\mytable{
{\tiny
\begin{tabular}{|c||c|c|c|c|c|c|c|c|c|c|c|c|c|c|c|c|c|c|c|c|}
\hline
N & 2 & 3 & 4 & 5 & 6 & 7 & 8 & 9 & 10 & 11 & 12 & 13 & 14 & 15 & 16 & 17 & 18\\
\hline
\hline
4 & 429 & 361 & 289 & 262 & 204 & 181 & 137 & 120 & 107 & 94 & 77 & 52 & 50 & 37 & 38 & 25 &
31\\
\hline
5 & 436 & 378 & 327 & 269 & 215 & 185 & 173 & 120 & 119 & 106 & 75 & 67 & 56 &  & 44 &  & 28\\
\hline
6 & 446 &  & 324 & 282 & 243 &  & 183 & 154 & 115 & 107 & 88 & 68 & 65 & 59 & 36 & 47 &\\
\hline
7 & 453 & 400 & 330 & 287 &  & 208 & 176 & 142 & 126 & 97 & 74 & 69 & 50 &  & 35 & 39 & 33\\
\hline
8 & 440 & 396 & 275 & 230 & 198 & 149 & 137 & 116 & 103 & 63 & 57 & 51 & 39 & 34 & 37 & 25 & 23\\
\hline
9 & 463 & 404 & 334 & 276 & 208 & 180 & 148 & 121 & 86 & 70 & 73 & 41 & 44 & 29 & 29 & 17 & 15\\
\hline
10 & 461 & 383 & 328 & 274 & 221 & 165 & 156 & 113 & 83 & 71 & 60 & 46 & 42 & 30 & 26 & 10 & 17\\
\hline
14 & 460 & 377 & 295 & 244 &  & 140 & 108 & 79 & 54 & 41 & 33 & 19 & 14 & 14 & 8 & 9 & 8\\
\hline
17 & 453 & 365 & 293 & 221 & 167 & 118 & 89 & 56 & 56 & 33 & 16 & 17 & 10 & 4 & 2 & 4 &\\
\hline
20 & 455 & 373 & 305 & 226 & 153 &  & 73 & 43 & 36 & 21 & 11 &  & 8 &  & 3 & 3 & 2\\
\hline
\end{tabular}}
\caption{Number of successes for finding $x$ out of $500$ tries}\label{find_x_10}
}
The results suggest that while we already obtain solutions (with nontrivial probability)
for some nontrivial parameters, we must extend the approach in order to consider harder
parameters. A successful extension is discussed in \cite{BraidEqns}. In the sequel
we discuss some other possible extensions.

\section{Possible improvements and conclusions}\label{future}
One approach is to create new conjugates by multiplying any
number of existing ones (or their inverses).
In fact, if $\cB$ is the group generated by $b_1,\dots,b_n$, then
the group generated by $xb_1x\inv,\dots,xb_nx\inv$ is
$x\cB x\inv$.
By Section \ref{many}, this does not help much.

The algorithm can be randomized by conjugating the given
elements $xb_1x\inv,\dots,xb_nx\inv$ by a random (known) element
$y\in \<a_1,\dots,a_m\>$,
so that running it several times increases the success probability.
The problem with this approach is that the conjugator becomes longer
and therefore the probability of success in each single case decreases.

Our experiments showed that the peeling off process often enters a loop,
that is, a stage to which we return every several steps.
This can sometimes be solved by conjugating with a random known element
after we enter the loop.
We also tried to change the length function or the linear ordering when we enter
a loop.
These approaches were successful for small parameters but
did not result in a significant improvement for large
parameters.

We did not try approaches of learning algorithms, neural networks, etc.
A simple example is to try and learn the distribution of the lengths for the
correct generator and define the linear ordering according to the likelihood
test.

The purpose of this paper was to check the applicability of
Hughes and Tannenbaum's length-based approach against the key agreement protocols
introduced in \cite{Anshel, Ko}. Our results suggest that this approach
requires an unfeasible computational power in order to solve the generalized conjugacy
search problem for the parameters used in these protocols.
However, this method has natural extensions which can make it applicable:
In \cite{BraidEqns} we suggest one particularly successful extension,
and it turns out that it can solve these and other problems
with standard computational power.


\begin{thebibliography}{10}

\bibitem{Anshel}
I.\ Anshel, M.\ Anshel and D.\ Goldfeld,
\emph{An algebraic method for public-key cryptography},
Mathematical Research Letters \textbf{6} (1999),
287--291.

\bibitem{Bangert}
P.\ D.\ Bangert,
\emph{Raid Braid: Fast Conjugacy Disassembly in Braid and Other Groups},
Applications of Computer Algebra (ACA-2004).

\bibitem{Cha}
J.\ C.\ Cha, K.\ H.\ Ko, S.\ J.\ Lee, J.\ W.\ Han and J.\ H.\ Cheon,
\emph{An Efficient Implementation of Braid Groups},
ASIACRYPT 2001, Lecture Notes in Computer Science \textbf{2248} (2001),
144--156.

\bibitem{Deh2}
P.\ Dehornoy,
\emph{A fast method for comparing braids},
Advances in Mathematics \textbf{125} (1997),
200--235.

\bibitem{DeBBC}
P.\ Dehornoy,
\emph{Braid-based cryptography},
in: \textbf{Group Theory, Statistics, and Cryptography},
Contemporary Mathematics \textbf{360} (2004),
5--33.

\bibitem{DehShift}
P.\ Dehornoy,
\emph{Using shifted conjugacy in braid-based cryptography},
Contemporary Mathematics \textbf{418} (2006), 65--73.

\bibitem{WPG}
D.\ B.\ A.\ Epstein, M.\ S.\ Paterson, G.\ W.\ Camon, D.\ F.\ Holt, S.\ V.\ Levy, and W.\ P.\ Thurston,
\textbf{Word Processing in Groups}, Jones and Bartlett, Boston, MA, 1992.

\bibitem{draft}
D.\ Garber, S.\ Kaplan, M.\ Teicher, B.\ Tsaban, and U.\ Vishne,
\emph{Length-based conjugacy search in the braid group}, Version 1,
Mathematics ArXiv eprint (2002)
\verb|math.GR/0209267 (v.1)|

\bibitem{BraidEqns}
D.\ Garber, S.\ Kaplan, M.\ Teicher, B.\ Tsaban, and U.\ Vishne,
\emph{Probabilistic solutions of equations in the braid group},
Advances in Applied Mathematics \textbf{35} (2005),
323--334.

\bibitem{Garrett}
P.\ Garrett,
\emph{Making, Braking Codes: Introduction to Cryptology},
Prentice Hall, 2000.

\bibitem{Gon}
 J.\ Goz\'alez-Meneses,
 \emph{Improving an algorithm to solve multiple simultaneous conjugacy problems in braid groups},
in: \textbf{Geometric methods in group theory}, Contemporary Mathematics \textbf{372} (2005),
35--42.

\bibitem{BGC}
Helger Lipmaa,
\emph{Cryptography and Braid Groups homepage},\\
\verb|http://www.tcs.hut.fi/~helger/crypto/link/public/braid/|

\bibitem{HofSte}
D.\ Hofheinz and R.\ Steinwandt,
\emph{A Practical Attack on Some Braid Group Based Cryptographic Primitives},
Proceedings of PKC 2003, Lecture Notes in Computer Science \textbf{2567} (2003),
187--198.

\bibitem{HT}
J.\ Hughes and A.\ Tannenbaum,
\emph{Length-based attacks for certain group based encryption rewriting systems},
Workshop SECI02 S\'ecurit\'e de la Communication sur Internet,
September 2002.

\bibitem{Ko}
K.\ H.\ Ko, S.\ J.\ Lee, J.\ H.\ Cheon, J.\ W.\ Han, S.\ J.\ Kang and C.\ S.\ Park,
\emph{New Public-key Cryptosystem using Braid Groups},
CRYPTO 2000, Lecture Notes in Computer Science \textbf{1880} (2000),
166--183.

\bibitem{Lee}
E.\ Lee,
\emph{Braid groups in cryptology},
IEICE Transactions on Fundamentals of Electronics, Communications and Computer Sciences \textbf{E87-A} (2004),
986--992.


\bibitem{LeePark}
E.\ Lee, J.\ H.\ Park,
\emph{Cryptanalysis of the public-key encryption based on braid groups},
EUROCRYPT 2003, Lecture Notes in Computer Science \textbf{2656} (2003),
477--490.

\bibitem{Paeng}
S.\ H.\ Paeng, K.\ C.\ Ha, J.\ H.\ Kim, S.\ Chee and C.\ Park,
\emph{New Public Key Cryptosystem Using Finite Non Abelian Groups},
CRYPTO 2001, Lecture Notes in Computer Science \textbf{2139} (2001),
470--485.

\bibitem{Shp}
V.\ Shpilrain,
\emph{Assessing security of some group based cryptosystems},
Contemporary Mathematics \textbf{360} (2004),
167--177.

\bibitem{UshShp}
A.\ Ushakov and V.\ Shpilrain,
\emph{The conjugacy search problem in public key cryptography: unnecessary and insufficient},
Applicable Algebra in Engineering, Communication and Computing \textbf{17} (2006), 285--289.

\bibitem{ShpZap}
V.\ Shpilrain and G.\ Zapata,
\emph{Combinatorial group theory and public key cryptography},
Applicable Algebra in Engineering, Communication and Computing \textbf{17} (2006), 291--302.

\end{thebibliography}
\end{document}